\documentclass[a4,psfig,english,11pt,epsf,portrait]{article}
\usepackage{amsmath,amsfonts,latexsym,amscd,hyperref}
\usepackage[english]{babel}
\usepackage{epsfig}
\usepackage[T1]{fontenc}
\usepackage[ansinew]{inputenc}  
\usepackage{graphics}
\usepackage{graphicx}
\usepackage{tikz}
\usepackage{amsmath}
\usepackage{amsthm}
\usepackage{amssymb}
\usepackage{enumerate}
\usepackage{multicol}
\usepackage{hyperref}

\newcommand{\abs}[1]{\lvert#1\rvert}

\newtheorem{theorem}{Theorem}[section]
\newtheorem{remark}{Remark}[section]
\newtheorem{proposition}[theorem]{Proposition}
\newtheorem{definition}{Definition}[section]
\newtheorem{lemma}[theorem]{Lemma}

\numberwithin{equation}{section}

\usepackage{color}

\def\N{\mathbb{N}}

\title{The Fujita exponent for a heat equation with mixed local and nonlocal nonlinearities on the Heisenberg group}

\author{
Zineb Sabbagh\\
{\it \small 
Department of Mathematics, University M'hamed Bougara Boumerdes, Algeria}\\
Ahmad Z. Fino\\
{\it \small 
College of Engineering and Technology, American University of the Middle East, Egaila 54200, Kuwait}\\
Mokhtar Kirane\\
{\it \small 
Department of Mathematics, College of Computing and Mathematical Sciences},\\
{\it \small Khalifa University, P.O. Box: 127788,  Abu Dhabi, UAE}\\
}

\date{}

\begin{document}
	\maketitle

\begin{abstract}
This article deals with the problems of local and global solvability for a semilinear heat equation on the Heisenberg group involving a mixed local and nonlocal nonlinearity. The characteristic features of such equations, arising from the interplay between the geometric structure of the Heisenberg group and the combined nonlinearity, are analyzed in detail. The need to distinguish between subcritical and supercritical regimes is identified and justified through rigorous analysis. On the basis of the study, the author suggests precise conditions under which local-in-time mild solutions exist uniquely for regular, nonnegative initial data. It is proved that global existence holds under appropriate growth restrictions on the nonlinear terms. To complement these results, it is shown, by employing the capacity method, that solutions cannot exist globally in time when the nonlinearity exceeds a critical threshold. As a result, the Fujita exponent is formulated and identified as the dividing line between global existence and finite-time blow-up. In addition, lifespan estimates were obtained in the supercritical regime, providing insight into how the size of the initial data influences the time of blow-up.
\end{abstract}

\medskip




\section{Introduction}
In this paper, we examines the global existence/nonexistence of mild solutions of the following semilinear heat equation on the Heisenberg group
\begin{equation}\label{1}
	\left\{\begin{array}{ll}
		\,\, \displaystyle {\partial_t u-\Delta_{\mathbb{H}}
			u =  \int_0^t(t-s)^{-\gamma}|u|^{p_1-1}u(s)\,\mathrm{d}s} + |u|^{p_2-1}u, &\displaystyle {\eta\in {\mathbb{H}}^{n},\,t>0,}\\
		{}\\
		\displaystyle{u(\eta,0)=  u_0(\eta),\qquad\qquad}&\displaystyle{\eta\in {\mathbb{H}}^{n},}
	\end{array}
	\right.
\end{equation}
where \( \Delta_{\mathbb{H}} \) denotes the sub-Laplacian on the Heisenberg group \( \mathbb{H}^n \), \( \gamma \in [0,1) \), and \( p_1, p_2 > 1 \). The initial data \( u_0 \in C_0(\mathbb{H}^n) \) is assumed to be nonnegative and continuous vanishing at infinity.

The study of semilinear heat equations and their qualitative properties such as global existence, blow-up, and lifespan of solutions has attracted substantial attention in the past few decades due to its fundamental role in mathematical physics, geometry, and biological modeling. A prototypical model is the classical semilinear heat equation
\begin{equation}\label{3}
 u_t=\Delta u +u^p,\qquad \text{in}\,\,\mathbb{R}^n.
 \end{equation}
for which Fujita~\cite{Fujita} identified the critical exponent \( p_F = 1 + 2/n \), below which all nontrivial nonnegative solutions blow up in finite time, and above which global existence may hold for small data. Since then, many authors have extended Fujita-type results to equations involving different types of nonlinearities and more complex geometric settings, see e.g. \cite{CDW,Souplet2001,Zhang1}. The Heisenberg group provides a natural setting for subelliptic diffusion, and the operator \( \Delta_{\mathbb{H}} \) exhibits strong hypoellipticity and sub-Riemannian scaling properties. In this framework,  Zhang \cite{Zhang}, Pohozaev-V\'eron \cite{Pohozaev}, and Pascucci \cite{Pascucci} derived critical Fujita-type exponents for the pure power-type problem
\begin{equation}\label{5}
\partial_t u - \Delta_{\mathbb{H}} u = |u|^{p-1} u,\qquad \text{in}\,\,\mathbb{H}^n.
 \end{equation}
with the critical threshold \( p_F = 1 + 2/Q\), where \( Q = 2n + 2 \) is the homogeneous dimension of \( \mathbb{H}^n \). Their works established the fundamental dichotomy between global existence and finite-time blow-up, similar to the Euclidean case.

To our knowledge, there were only a few papers addressing the study of diffusion equations with mixed local-nonlocal nonlinearities. The inclusion of both a nonlocal memory term involving a time-convolution and a power-type reaction term makes problem \eqref{1} significantly richer and more challenging. This type of mixed local-nonlocal nonlinearity has been explored in various settings. For instance, in the Euclidean setting, Zhang \cite{Zhang1}  considered the following Caputo time-fractional diffusion equation with memory and reaction terms in a bounded domain
$$
\begin{array}{ll}
\displaystyle {}^C_0D^\alpha_t u-\Delta u= \int_0^t(t-s)^{-\gamma}|u|^{p-1}u(s)\,\mathrm{d}s - |u|^{q-1}u,&\qquad {x\in \Omega\subset \mathbb{R}^n,\,\,\,t>0,}\end{array}
$$
where \( \alpha \in (0,1] \), $\gamma\in[0,1)$, $p,q>1$,  and studied blow-up and global existence behavior depending on the relation between \( p \) and \( q \). When $\alpha = 1$ and $\gamma =0$,  Souplet~\cite{Souplet2001} proved that if $p > q$ and $u_0 \geq 0$, then all nontrivial solutions blow up in finite time in the $L^{\infty}$ norm, while if $p \leq q$, and $u_0 \geq 0$, then the nontrivial solution is global and unbounded. When \( \alpha = 1\) and the reaction term is omitted, the problem in the Euclidean space reduces to the purely nonlocal equation
\begin{equation}\label{2}
\begin{array}{ll}
\displaystyle \partial_t u-\Delta u= \int_0^t(t-s)^{-\gamma}|u|^{p-1}u(s)\,\mathrm{d}s,&\qquad {x\in \mathbb{R}^n,\,\,\,t>0.}\end{array}
\end{equation} 
In this setting, Cazenave et al.~\cite{CDW} showed that the critical exponent governing the global behavior of solutions is given by 
$$p_{\mathrm{crit}}=\max\Big\{\frac{1}{\gamma},1+\frac{2(2-\gamma)}{(n-2+2\gamma)_+}\Big\}\in(0,+\infty].$$
Specifically, it was shown that if $p \leq p_{\mathrm{crit}}$, then all nontrivial nonnegative solutions blow up in finite time, while for $p>p_{\mathrm{crit}}$, there exists a unique global solution for sufficiently small initial data.

Analogous results were obtained in the non-Euclidean setting of the Heisenberg group. Recently, Fino et al. \cite{FinoKiraneTorebek} considered the equation
 \begin{equation}\label{4}
\begin{array}{ll}
\displaystyle \partial_t u-\Delta u= \int_0^t(t-s)^{-\gamma}|u|^{p-1}u(s)\,\mathrm{d}s,&\qquad {\eta\in \mathbb{H}^n,\,\,t>0,}\end{array}
\end{equation} 
 and proved that the Fujita-type critical exponent is
 $$p_{\mathrm{crit}}=\max\Big\{\frac{1}{\gamma},1+\frac{2(2-\gamma)}{Q-2+2\gamma}\Big\}\in(0,+\infty].$$
Moreover, as $\gamma\to1$, using the relation
$$\lim_{\gamma\to1}c_\gamma\,s_+^{-\gamma}=\delta_0(s)\quad\text{in distributional sense with}\,\,s_+^{-\gamma}:=\left\{\begin{array}{ll}
s^{-\gamma}&\,\,\text{if}\,\,s>0,\\
0&\,\,\text{if}\,\,s<0,
\end{array}
\right.
$$
where $c_\gamma=1/\Gamma(1-\gamma)$, and through a suitable change of variables, we observe that equation \eqref{2} (in $\mathbb{R}^n$) and equation \eqref{4} (in $\mathbb{H}^n$)  reduce to the classical semilinear heat equation~\eqref{3} and \eqref{5} respectively, thus recovering Fujita's result.

Further extensions of Fujita-type results on the Heisenberg group have been obtained by many authors (see for example \cite{DAmbrosio,FRT,Georgiev,Kirane3, Kirane5} and the references therein). 

Motivated by these observations, we emphasize the relevance of Fujita exponents for heat equations \eqref{1} involving the Heisenberg operator with mixed local and nonlocal nonlinearities

For the convenience of the readers, we present the main results of the paper in the subsection below, and then in the following sections we present the proofs of the main results in turn.

\subsection*{Main results}
\begin{definition}[Mild solution]${}$\\
Let $u_0\in C_0(\mathbb{H}^n)$, $n\geq1$, $p_i>1$ for $i=1,2$, and $T>0$. We say that $u\in C([0,T],C_0(\mathbb{H}^n))$
is a mild solution of problem \eqref{1} if $u$ satisfies the following integral equation
\begin{equation}\label{IE}
    u(t,\eta)=S_{\mathbb{H}}(t) u_0(\eta)+\Gamma(\alpha)\int_{0}^tS_{\mathbb{H}}(t-s)I_{0|s}^\alpha(|u|^{p_1-1}u)(s,\eta)\,ds+\int_{0}^tS_{\mathbb{H}}(t-s)|u|^{p_2-1}u(s,\eta)\,ds,
\end{equation}
for all $\eta\in \mathbb{H}^n,\,t\in[0,T]$, where $\alpha=1-\gamma$, $\Gamma$ is the Euler gamma function and $I_{0|s}^\alpha$ is the Riemann-Liouville left-sided fractional integral defined in \eqref{I1} below. In addition, if $u$ is a mild solution of \eqref{1} in $[0,T]$ for all $T>0$, then $u$ is called global-in-time mild solution of \eqref{1}.
\end{definition}
First, by applying the same strategies as in \cite[Theorem~1]{FinoKiraneTorebek}, we establish the existence and uniqueness of a local mild solution.
\begin{theorem}[Local existence]\label{Local} Let $u_0\in C_0(\mathbb{H}^n)$, with $n\geq1$ and $p_i>1$ for $i=1,2$. Then there exists a maximal
time $T_{\max}>0$ and a unique mild solution $u\in
C([0,T_{\max}),C_0(\mathbb{H}^n))$ to problem \eqref{1}. Moreover, we have the alternative:\\
\begin{itemize} 
\item[$\bullet$] Either $T_{\max}=+\infty$, in which case the solution is global, or
\item[$\bullet$]  $T_{\max}<+\infty$ and 
$$\liminf_{t\rightarrow T_{\max}}\|u\|_{L^\infty((0,t),L^\infty(\mathbb{H}^{n}))}=+\infty,$$
that is, the solution blows up in finite time.
\end{itemize} 
 In addition, if $u_0\geq0,$ $u_0\not\equiv0,$ then the solution remains strictly positive for all $0<t<T_{\max}.$
 Furthermore, if $u_0\in L^r(\mathbb{H}^n),$ for some $1\leq r<\infty,$ then $u\in C([0,T_{\max}),L^r(\mathbb{H}^n))$.
\end{theorem}
Next, we identify conditions that guarantee the global existence of mild solutions and, conversely, characterize the parameter ranges leading to finite-time blow-up. Before stating the theorems, recall that $Q=2n+2$ denotes the homogeneous dimension of the Heisenberg group. We now introduce the critical exponents as follows:
\begin{equation}\label{critexp}
p_1^* = \max\left\{\frac{1}{\gamma}, p_\gamma\right\}, \qquad p_2^* = 1+ \frac{2}{Q} ,\qquad\hbox{and}\qquad p_2^{**}=\max\left\{\frac{\gamma-\gamma^2+1}{\gamma(2-\gamma)} , 1+\frac{2}{Q-2+2\gamma}\right\} ,
\end{equation}
where
\begin{eqnarray*}
	p_\gamma =	1+ \frac{2(2-\gamma)}{Q-2+2\gamma} .
\end{eqnarray*}
We also define $\tilde{p}_2:= (p_1+1-\gamma)/(2-\gamma)$ and the scaled exponents as
$$q_{\mathrm{sc}}=\left\{\begin{array}{ll}
 \displaystyle \frac{Q(p_1-1)}{2(2-\gamma)},&\qquad\hbox{if}\,\, p_2\geq \tilde{p}_2,\\\\
\displaystyle \frac{Q(p_2-1)}{2},&\qquad\hbox{if}\,\, p_2\leq \tilde{p}_2.\\
\end{array}
\right.$$
In addition, in the case when $p_2>p_2^*$, the strict inequality  $p_2^*<p_2^{**}$ allows us to distinguish between two subcases: $p_2>p_2^{**}$ and $p_2^*<p_2\leq p_2^{**}$.
\begin{theorem}[Global existence]\label{global1}${}$\\ 
Let $\gamma\in(0,1)$, $p_i>1$ for $i=1,2$, and let $u_0\in  C_0(\mathbb{H}^{n})\cap L^{q_{\mathrm{sc}}}(\mathbb{H}^{n})$ with $\|u_0\|_{L^\infty}$ and $\|u_0\|_{ L^{q_{\mathrm{sc}}}}$ sufficiently small. If 
$$p_1>p_1^*\qquad\hbox{and}\qquad p_2>p_2^{**},$$
 then problem \eqref{1} admits a unique global mild solution 
$$u\in C([0,\infty);C_0(\mathbb{H}^{n})\cap L^{q_{\mathrm{sc}}}(\mathbb{H}^{n})).$$ 
\end{theorem}
\begin{remark}${}$
\begin{itemize}
\item[a)] The condition $p_1>p_1^*$ is equivalent to $\tilde{p}_2>p_2^{**}$, which naturally leads to considering three distinct subcases in the proof of Theorem \ref{global1}: $p_2=\tilde{p}_2$, $p_2>\tilde{p}_2$, and $p_2^{**}<p_2<\tilde{p}_2$. In the case $p_2= \tilde{p}_2$, it suffices to assume that $\|u_0\|_{ L^{q_{\mathrm{sc}}}}$ is sufficiently small.
\item[b)] Observe that, when $p_2\geq \tilde{p}_2$, one may assume the decay condition $|u_0(\eta)|\leq C
(1+|\eta|_{_\mathbb{H}})^{-\kappa}$ with $\kappa>2(2-\gamma)/(p_1-1)$ instead of the integrability condition $u_0\in L^{q_{\mathrm{sc}}}(\mathbb{H}^n)$. This example highlights the relationship between the decay rate of $u_0$ and that described in Theorem \ref{Lifespan-2} below.
\item[c)] We note that Theorems $\ref{Local}$, and $\ref{global1}$ remain valid when the nonlinearities of the form $u\mapsto |u|^{p_i-1}u$ (for $i=1,2$) are replaced by more general nonlinear functions $f$ and $g$ satisfying the growth conditions given in \eqref{growth1} below.
\end{itemize}
\end{remark}

\begin{theorem}[Blow-up]\label{T4} Let $0\leq \gamma<1$, $u_0\in C_0(\mathbb{H}^n)$ be such that $u_0\geq0$ and $u_0\not\equiv0$. If
	$$p_1\leq p_1^*\qquad\hbox{or}\qquad p_2\leq p_2^*,$$
	then the mild solution of  \eqref{1} blows-up in finite time.
\end{theorem}
\begin{remark} The question of whether global existence or blow-up occurs remains open in the case where $p_1>p_1^*$ and $p_2^*<p_2\leq  p_2^{**}$.
\end{remark} 

Next, we investigate estimates for the lifespan of local-in-time weak or mild solutions to equation \eqref{1}. Specifically, we aim to understand how the maximal existence time depends on the size of the initial data. To capture this dependency, we consider scaled initial data of the form \( u_0^\varepsilon = \varepsilon u_0 \), where \( \varepsilon > 0 \) is a small parameter and \textcolor{red}{$u_0\in C_0(\mathbb{H}^n)\cap L^1(\mathbb{H}^n)$.} We define the lifespan \( T_\varepsilon \) as the supremum of all \( T > 0 \) such that the corresponding solution \( u^\varepsilon(t, \eta) \) exists and remains in the appropriate function space on the time interval \( [0, T) \). Our goal is to derive asymptotic estimates on \( T_\varepsilon \) as \( \varepsilon \to 0^+ \), thereby revealing how the smallness of the initial data governs the time interval of existence.

The proofs of the following theorems follow directly from the arguments presented in \cite[Theorems~3-4]{FinoKiraneTorebek} and \cite[Theorem~1]{Georgiev}
\begin{theorem}[Lifespan estimates]\label{Lifespan} Let $0\leq\gamma<1$, $p_i>1$ for $i=1,2$, $u_0\in C_0(\mathbb{H}^n)\cap L^1(\mathbb{H})$, and $u_0\geq0$ such that $\displaystyle \int_{\mathbb{H}^n} u_0(\eta)\,d\eta>0$. Then there exists $\varepsilon_0>0$ such that for any $\varepsilon\in(0,\varepsilon_0]$ it holds
$$T_\varepsilon\leq\left\{\begin{array}{ll}
C \varepsilon^{-\left(\frac{2-\gamma}{p_1-1}-\frac{Q}{2}\right)^{-1}}&\quad\hbox{if}\quad p_1<p_{\mathrm{sc}_1}:=1+\frac{2(2-\gamma)}{Q},\\
 C \varepsilon^{-\left(\frac{1}{p_2-1}-\frac{Q}{2}\right)^{-1}}&\quad\hbox{if}\quad p_2<p^*_2,\\
 C \exp(\varepsilon^{-(p_2-1)})&\quad\hbox{if}\quad p_2=p^*_2,\\
 \end{array}
\right.
 $$
where $C>0$ is a positive constant independent of $\varepsilon$.
\end{theorem}
The result on the lifespan can be relaxed in the following manner.
\begin{theorem}[Lifespan estimates]\label{Lifespan-2} Let $0\leq\gamma<1$, $p_i>1$ for $i=1,2$, and $u_0\in C_0(\mathbb{H}^n)\cap L^1_{\mathrm{loc}}(\mathbb{H})$ such that 
\begin{equation}\label{corA}
u_0(\eta)\geq (1+|\eta|)^{-\kappa},\qquad\hbox{for all}\,\,\eta\in\mathbb{H}^n.
\end{equation}
 If $0<\kappa<2(2-\gamma)/(p_1-1)$,  there exists a constant $\varepsilon_0>0$ such that the lifespan $T_\varepsilon$ verifies
$$
T_\varepsilon \leq\,C \, \varepsilon^{-\left(\frac{2-\gamma}{p_1-1}-\frac{\kappa}{2}\right)^{-1}}, \qquad \hbox{for any}\,\, \varepsilon \in ( 0,\varepsilon_0].
$$
\end{theorem}
\begin{remark} Theorem \ref{Lifespan} does not cover lifespan estimates for all blow-up scenarios. In particular, for the range $p_{\mathrm{sc}_1}\leq p_1\leq p^*_1$, the lifespan remains unknown. However, Theorem \ref{Lifespan-2} provides a relaxed estimate for the lifespan in some of these cases, partially addressing this gap.
\end{remark}
The proofs of local and global existence of solutions rely on $L^p-L^q$ estimates and the application of a fixed point theorem. In contrast, the analysis of blow-up and lifespan of solutions is carried out using nonlinear capacity estimates, also known as the rescaled test function method (cf. \cite{Fino3,Fino5,Kirane5,19,Zhang}). 

The paper is organized as follows. Section \ref{sec2} introduces key definitions, terminology, and preliminary results related to the Heisenberg group, the sub-Laplacian operator, and its heat kernel. Sections \ref{sec3} and \ref{sec4} are dedicated to proving the global existence and blow-up results.


\section{Preliminaries}\label{sec2}
As a preliminary step, we recall the definition and fundamental properties of the Heisenberg group (see, e.g., \cite{Folland,Follandstein}).

\subsection{Heisenberg group}
The Heisenberg group $\mathbb{H}^n$  is the space $\mathbb{R}^{2n+1}=\mathbb{R}^n\times \mathbb{R}^n\times \mathbb{R}$ equipped with the group operation
$$\eta\circ\eta^\prime=(x+x^\prime,y+y^\prime,\tau+\tau^\prime+2(x\cdotp y^\prime-x^\prime\cdotp y)),$$
where $\eta=(x,y,\tau)$, $\eta^\prime=(x^\prime,y^\prime,\tau^\prime)$, and $\cdotp$ is the standard scalar product in $\mathbb{R}^n$. Let us denote the parabolic dilation in $\mathbb{R}^{2n+1}$ by $\delta_\lambda$, namely, $\delta_\lambda(\eta)=(\lambda x, \lambda y, \lambda^2 \tau)$ for any $\lambda>0$, $\eta=(x,y,\tau)\in\mathbb{H}^n$. The Jacobian determinant of $\delta_\lambda$ is $\lambda^Q$, where $Q=2n+2$ is the homogeneous dimension of $\mathbb{H}^n$. A direct calculation shows that $\delta_\lambda$ is an automorphism of $\mathbb{H}^n$ for every $\lambda>0$, and therefore  $\mathbb{H}^n=(\mathbb{R}^{2n+1},\circ,\delta_\lambda)$ is a homogeneous Lie group on $\mathbb{R}^{2n+1}$.\\
The homogeneous Heisenberg norm (also called Kor\'anyi norm) is derived from an anisotropic dilation on the Heisenberg group and is defined by
$$|\eta|_{_{\mathbb{H}}}=\left(\left(\sum_{i=1}^n (x_i^2+y_i^2)\right)^2+\tau^2\right)^{\frac{1}{4}}=\left((|x|^2+|y|^2)^2+\tau^2\right)^{\frac{1}{4}},$$
where $|\cdotp|$ is the Euclidean norm associated with $\mathbb{R}^n$. The associated Kor\'anyi distance between two points $\eta$ and $\xi$ of $\mathbb{H}$ is defined by $d_{_{\mathbb{H}}}(\eta,\xi)=|\xi^{-1}\circ\eta|_{_{\mathbb{H}}},\quad \eta,\xi\in\mathbb{H}$, where $\xi^{-1}$ denotes the inverse of $\xi$ with respect to the group action, i.e. $\xi^{-1}=-\xi$. This metric induces a topology on $\mathbb{H}^n$.\\
The Heisenberg convolution between any two regular functions $f$ and $g$ is defined by
$$
(f\ast_{_{\mathbb{H}}}g)(\eta)=\int_{\mathbb{H}^n}f(\eta\circ\xi^{-1})g(\xi)\,d\xi=\int_{\mathbb{H}^n}f(\xi)g(\xi^{-1}\circ\eta)\,d\xi.$$
The left-invariant vector fields that span the Lie algebra \(\mathfrak{h}_n\) of \(\mathbb{H}^n\) are given by
\[
X_i = \partial_{x_i}-2 y_i\partial_\tau, \quad
Y_i = \partial_{y_i}+2 x_i\partial_\tau, \quad
T = \partial_\tau, \quad i = 1, \ldots, n.
\]
The Heisenberg gradient is given by
\begin{equation}\label{48}
\nabla_{\mathbb{H}}=(X_1,\dots,X_n,Y_1,\dots,Y_n),
\end{equation}
and the sub-Laplacian (also referred to as Kohn Laplacian) is defined as
\begin{equation}\label{40}
\Delta_{\mathbb{H}}=\sum_{i=1}^n(X_i^2+Y_i^2)=\Delta_x+\Delta_y+4(|x|^2+|y|^2)\partial_\tau^2+4\sum_{i=1}^{n}\left(x_i\partial_{y_i\tau}^2-y_i\partial_{x_i\tau}^2\right),
\end{equation}
where $\Delta_x=\nabla_x\cdotp\nabla_x$ and $\Delta_y=\nabla_y\cdotp\nabla_y$ stand for the Laplace operators on $\mathbb{R}^n$.

\subsection{Heat kernel}
We recall the definition and some properties related to the heat kernel associated with $-\Delta_{\mathbb{H}}$ on the Heisenberg group.
\begin{proposition}\label{Heat}\cite[Theorem~1.8]{BahouriGallagher}${}$\\
There exists a function $h\in\mathcal{S}(\mathbb{H}^n)$ such that if u denotes the
solution of the free heat equation on the Heisenberg group
\begin{equation}\label{}
\begin{array}{ll}
\partial_t u-\Delta_{\mathbb{H}}u=0,&\qquad {\eta\in \mathbb{H}^n,\,\,t>0,}
 \\{}\\ u(\eta,0)=u_{0}(\eta), &\qquad \eta\in \mathbb{H}^n,
 \end{array}
\end{equation}
then we have
$$u(\cdotp,t)=h_t\ast_{_{\mathbb{H}}} u_0,$$
where $\ast_{_{\mathbb{H}}}$ denotes the convolution on the Heisenberg group defined above, while the heat kernel $h_t$ associated with $-\Delta_{\mathbb{H}}$ is defined by
$$h_t(\eta)=\frac{1}{t^{n+1}}h\left(\frac{x}{\sqrt{t}},\frac{y}{\sqrt{t}},\frac{s}{t} \right),\qquad\hbox{for all}\,\,\,\eta=(x,y,s)\in\mathbb{H}^n, \,\,t>0.$$
\end{proposition} 
The following proposition follows directly from \cite[Theorem~3.1]{Folland}, \cite[Corollary~1.70]{Follandstein}, \cite[Proposition~1.68]{Follandstein}, and \cite[Corollary~2.3]{Pazy}.
\begin{proposition}\label{properties}${}$\\
There is a unique semigroup $(S_{\mathbb{H}}(t))_{t>0}$ generated by $\Delta_{\mathbb{H}}$ and satisfying the following properties:
\begin{itemize}
\item[$\mathrm{(i)}$] $(S_{\mathbb{H}}(t))_{t>0}$ is a contraction semigroup on $L^p(\mathbb{H}^n)$, $1\leq p\leq \infty$, which is strongly continuous for $p<\infty$;
\item[$\mathrm{(ii)}$] $(S_{\mathbb{H}}(t))_{t>0}$  is a strongly continuous semigroup on $C_0(\mathbb{H}^n)$;
\item[$\mathrm{(iii)}$] For every $v\in X$, where $X$ is either $L^p(\mathbb{H}^n)$ for $1\leq p< \infty$ or $C_0(\mathbb{H}^n)$, the map $t\longmapsto S_{\mathbb{H}}(t)v$ is continuous from $[0,\infty)$ into $X$;
\item[$\mathrm{(iv)}$] $S_{\mathbb{H}}(t)f=h_t\ast_{_{\mathbb{H}}} f$, for all $f\in L^p(\mathbb{H}^n)$, $1\leq p\leq \infty$, $t>0$;
\item[$\mathrm{(v)}$] $\displaystyle\int_{\mathbb{H}^n}h_t(\eta)\,d\eta=1$, for all $t>0$;
\item[$\mathrm{(vi)}$] $S_\mathbb{H}(t)\varphi\geq 0$ whenever $\varphi\geq0$, for all $t\geq0$. More precisely, the heat kernel satisfies
$$h_t(\eta)>0,\qquad \hbox{for all}\,\, \eta \in \mathbb{H}^n, t>0;$$ 
\item[$\mathrm{(vii)}$] $h_{r^2t}(rx,ry,r^2s)=r^{-Q}h_{t}(x,y,s)$, for all $r,t>0$, $(x, y,s) \in \mathbb{H}^n$; 
\item[$\mathrm{(viii)}$] $h_t(\eta)=h_t(\eta^{-1})$, for all $\eta \in \mathbb{H}^n$, $t>0$;
\item[$\mathrm{(ix)}$] $h_t\ast_{_{\mathbb{H}}} h_s=h_{t+s}$, for all $s,t>0$. 
\end{itemize}
\end{proposition} 
\begin{lemma}\label{Lp-Lqestimate}($L^p-L^q$ estimate)\cite[Lemma~3]{FinoKiraneTorebek}${}$\\
 Suppose $1\leq p\leq q\leq \infty$. Then there exists a positive constant $C$ such that for every $f\in L^p(\mathbb{H}^n)$, the following inequality holds:
\begin{equation}\label{semigroup}
\|S_{\mathbb{H}}(t)f\|_{L^q}\leq C\,t^{-\frac{Q}{2}(\frac{1}{p}-\frac{1}{q})}\|f\|_{L^p},\qquad t>0.
\end{equation}
\end{lemma}
\begin{lemma}\label{lemmaSg}
	Let $1\leq q<r<\infty$, and $\alpha = \frac{Q}{2}\left(\frac{1}{q} - \frac{1}{r}\right)$. Then
	\begin{equation}
		\lim\limits_{t\rightarrow 0^+}t^{\alpha}\|S_\mathbb{H}(t)\varphi\|_{L^r}=0,
	\end{equation}
	for any $\varphi\in L^q(\mathbb{H}^{n})$.
\end{lemma}
\proof Given $\varphi \in L^q(\mathbb{H}^{n})$. Then, by density argument, there exists a sequence $(\varphi_n)_{n \in \mathbb{N}}\subset C_c^\infty(\mathbb{H}^{n})$ converging to $\varphi$ in $L^q(\mathbb{H}^{n})$. Using the semigroup decay \eqref{semigroup}, we have
	\begin{eqnarray*}
		t^\alpha\|S_\mathbb{H}(t)\varphi\|_{L^r}&\leq& t^\alpha\|S_\mathbb{H}(t)(\varphi-\varphi_n)\|_{L^r}+t^\alpha\|S_\mathbb{H}(t)\varphi_n\|_{L^r}\\
		&\leq&C\|\varphi-\varphi_n\|_{L^q}+t^\alpha\|\varphi_n\|_{L^r}\\
	\end{eqnarray*}
	for all $n\in \mathbb{N}, t>0$.
	Taking the limit as $t \rightarrow 0$ , we get
	\begin{equation*}
		\lim\limits_{t\rightarrow 0^+}t^{\alpha}\|S_\mathbb{H}(t)\varphi\|_{L^r}\leq C\|\varphi-\varphi_n\|_{L^q},
	\end{equation*}
	and the result follows by letting $n \rightarrow \infty$. \hfill$\square$\\
\begin{definition}
Let $f\in L^1(c,d)$, $-\infty<c<d<\infty$. The left-sided Riemann-Liouville fractional integral is defined by
\begin{equation}\label{I1}
I^\alpha_{c|t}f(t):=\frac{1}{\Gamma(\alpha)}\int_{c}^t(t-s)^{-(1-\alpha)}f(s)\,ds, \quad t>c, \; \alpha\in(0,1),
\end{equation}
where $\Gamma$ is the Euler gamma function.
\end{definition}

\subsection{Comparison Principle}\label{comparisonprinciple}

In this subsection, we establish a comparison principle for mild solutions of \eqref{1}, where the nonlinearities $u\mapsto |u|^{p_i-1}u$ (for $i=1,2$) are replaced by general functions $f$ and $g$ satisfying the following growth conditions:
\begin{equation}\label{growth1}
	|f(u)-f(v)| \leq c|u-v|(|u|^{p_1-1} + |v|^{p_1-1}),\quad\quad |g(u)-g(v)| \leq c|u-v|(|u|^{p_2-1} + |v|^{p_2-1}).
\end{equation}
Given an initial condition $u_0$, we denote by $u(t;u_0,f,g)$ the corresponding mild solution of \eqref{1}, emphasizing its dependence on the initial data and the nonlinear terms.

 We are now in a position to state a comparison principle, motivated by the Picard iteration scheme associated with mild solutions. To this end, we first introduce the notion of a mild supersolution, which will play a key role in the formulation of the principle. A function $v:[0,T]\to C_0(\mathbb{H}^n)$ is called a mild supersolution of equation \eqref{1}, if it satisfies $v(t) \geq \Phi_{fg}(v)(t)$, for all $t\in[0,T]$, where the operator $\Phi_{fg}$ is defined in equation \eqref{lambdafg} below.

\begin{theorem}[Comparison principle]\label{comparison}${}$\\
	Let $u_0,v_0\in C_0(\mathbb{H}^n)$, and let $f,g, \tilde{f}, \tilde{g}$ be functions satisfying the growth conditions in  \eqref{growth1}.
	\begin{itemize}
		\item[$\mathrm{(i)}$] If $0\leq u_0\leq v_0$, and $f(\mu)\leq \tilde{f}(\nu)$, $g(\mu)\leq \tilde{g}(\nu)$, for all $0\leq\mu\leq \nu$, then 
		$$0\leq u(t;u_0,f,g) \leq v(t;v_0,\tilde{f},\tilde{g}),\quad\text{for all}\,\, t\in[0,\tilde{T}],$$
		 where $\tilde{T}>0$ is a common time of existence for $u$ and $v$.
	\item[$\mathrm{(ii)}$]  If $v$ is a mild supersolution of \eqref{1}, and $u$ is a mild solution with $u_0\geq0$, then $0\leq u\leq v$ for all $t\in[0,T]$.
	\end{itemize}
\end{theorem}

\proof We begin by defining the Picard sequences:
$$u_{1}(t) =S_\mathbb{H}(t)u_0,\qquad u_{k}(t) =\Phi_{fg} (u_{k-1})(t),\,\,\, k\geq2 ,$$
$$v_{1}(t) =S_\mathbb{H}(t)v_0 ,\qquad v_{k}(t) =\widetilde{\Phi}_{\tilde{f}\tilde{g}} (v_{k-1})(t),\,\,\, k\geq2 , $$
where the operators $\Phi_{fg}$ and $\widetilde{\Phi}_{\tilde{f}\tilde{g}}$ are defined by
	\begin{equation}\label{lambdafg}
	\Phi_{fg}(u)(t):=S_{\mathbb{H}}(t)u_0 + \int_{0}^{t}S_{\mathbb{H}}(t-s) \int_0^s(s-\tau)^{-\gamma}f(u(\tau))\, \mathrm{d}\tau \,\mathrm{d}s + \int_{0}^{t}S_{\mathbb{H}}(t-s)g(u(s)) \mathrm{d}s,
	\end{equation}
	$$\widetilde{\Phi}_{\tilde{f}\tilde{g}}(u)(t):=S_{\mathbb{H}}(t)v_0 + \int_{0}^{t}S_{\mathbb{H}}(t-s) \int_0^s(s-\tau)^{-\gamma}\tilde{f}(u(\tau))\, \mathrm{d}\tau \,\mathrm{d}s + \int_{0}^{t}S_{\mathbb{H}}(t-s)\tilde{g}(u(s)) \mathrm{d}s.$$
	By Proposition \ref{properties}-(vi), the semigroup $S_\mathbb{H}(t)$ preserves order:
	\begin{equation}\label{monotonicity}
		S_\mathbb{H}( t)\varphi \leq S_\mathbb{H}(t)\psi,\quad \hbox{whenever}\,\,\, \varphi\leq \psi.
	\end{equation} 
Applying this property yields $0\leq u_1(t)\leq v_1(t)$. We now proceed by induction to prove tha $0\leq u_k(t)\leq v_k(t)$, for all $k\in\N$. Assume this holds for some fixed $k$. Then, using the monotonicity \eqref{monotonicity} of $S_\mathbb{H}(t)$  along with the growth assumptions on $f, \tilde{f},g,\tilde{g}$, we obtain
$$0\leq u_{k+1}(t)=\Phi_{fg}(u_k)(t)\leq \widetilde{\Phi}_{\tilde{f}\tilde{g}}(v_k)(t)=v_{k+1}(t).$$
Hence, the inequality is preserved at each iteration step. Moreover, since the Picard sequences of the mild solutions converge uniformly on compact intervals (and in particular, pointwise), we conclude that
	$$u(t;u_0,f,g) = \lim\limits_{n\rightarrow\infty} u_n(t)\leq \lim\limits_{n\rightarrow\infty} v_n(t)= v(t;v_0,\tilde{f},\tilde{g}).$$
This completes the proof of part $\mathrm{(i)}$. To prove part $\mathrm{(ii)}$, we proceed by constructing a non-increasing Picard sequence $(w_k)_k$ starting from the mild supersolution $v$. Define
$$
		w_{1}(t) =v(t),\qquad\quad w_{k}(t) =\Phi_{fg} (w_{k-1})(t),\quad k\geq2.
$$
	By part $\mathrm{(i)}$, we have
	 $$w_{2}(t) =\Phi_{fg} (w_1)(t)=\Phi_{fg} (v)(t)\leq v(t)=w_1(t).$$ 
	An induction argument, combined with the monotonicity of the operator $\Phi_{fg}$ (as established in part $\mathrm{(i)}$), shows that the sequence $(w_k)_k$  is non-increasing.
	Furthermore, since $f$ and $g$ are nonnegative and the semigroup $S_\mathbb{H}(t)$ preserves nonnegativity, each function $w_k(t)$ remains nonnegative and satisfies $0\leq w_k(t)\leq v(t)$ for all $k\geq1$. 
	As a consequence, $(w_k)_k$ is a bounded, non-increasing sequence in the Banach space $C([0,T],C_0(\mathbb{H}^{n})))$, and therefore converges uniformly on compact subsets of $[0,T]$ to a limit function $w$ within the same space. By standard monotone convergence theorems in Banach spaces, this convergence is strong. Passing to the limit in the iteration yields
	 $$w(t)= \Phi_{fg} (w)(t)\leq v(t),$$
which shows that $w$ is a mild solution of problem \eqref{1}, and moreover satisfies $w\leq v$. Finally, by the uniqueness of mild solutions, we conclude that $u=w\leq v$, which completes the proof of part $\mathrm{(ii)}$. \hfill$\square$

\section{Global existence}\label{sec3}
We now turn to the proof of the global existence result stated in Theorem \ref{global1}. Since $p_1>p_1^*$ is equivalent to $\tilde{p}_2>p_2^{**}$, the argument is divided into the three cases: $p_2=\tilde{p}_2$, $p_2>\tilde{p}_2$, and $p_2^{**}<p_2<\tilde{p}_2$.\\
\noindent {\it The Case $p_2= \tilde{p}_2$.} The proof is divided into two steps.\\
{\it Step 1. Existence.} As $p_1>p_1^*$, it is possible to choose a positive constant $q>0$ such that
\begin{equation}\label{estiA}
    \frac{2-\gamma}{p_1-1}-\frac{1}{p_1}<\frac{Q}{2
    q}<\frac{1}{p_1-1},\quad\hbox{with}\,\, q\geq p_1.
\end{equation}
From this inequality, it follows that
\begin{equation}\label{estiB}
    q>\frac{Q(p_1-1)}{2}>q_{sc}>1.
\end{equation}
Next, define the parameter $\beta$ by
\begin{equation}\label{estiC}
    \beta:=\frac{Q}{2 q_{sc}}-\frac{Q}{2
    q}=\frac{2-\gamma}{p_1-1}-\frac{Q}{2
    q}=\frac{1}{p_2-1}-\frac{Q}{2
    q}.
\end{equation}
Using the relations (\ref{estiA})-(\ref{estiC}) along with the condition $p_2<p_1$, one can verify that $p_i\beta<1$ for $i=1,2$, $q> p_2$, as well as the relations
\begin{equation}\label{estiD}
    \beta>\frac{1-\gamma}{p_1-1}>0,\quad \frac{Q(p_1-1)}{2
    q}+(p_1-1)\beta+\gamma=2,\quad \frac{Q(p_2-1)}{2
    q}+(p_2-1)\beta=1,
\end{equation}
hold. Since $u_0\in L^{q_{\mathrm{sc}}}(\mathbb{H}^{n})$, and using inequality \eqref{semigroup} along with the conditions $q>q_{\mathrm{sc}}$ and relation (\ref{estiC}), we
get
\begin{equation}\label{estiE}
    \sup_{t>0}t^\beta\|S_{\mathbb{H}}(t)u_0\|_{L^q}\leq
   C \|u_0\|_{L^{q_{\mathrm{sc}}}}=\eta<\infty,
\end{equation}
and
\begin{equation}\label{estiEE}
    \|S_{\mathbb{H}}(t)u_0\|_{L^{q_{\mathrm{sc}}}}\leq
   C \|u_0\|_{L^{q_{\mathrm{sc}}}}=\eta<\infty.
\end{equation}
We define the set
\begin{equation}\label{estiF}
    \Xi:=\left\{u\in
    L^\infty((0,\infty),L^q(\mathbb{H}^{n})\cap L^{q_{\mathrm{sc}}}(\mathbb{H}^{n}));\;\sup_{t>0}t^\beta\|u(t)\|_{L^q}\leq\delta\,\,\hbox{and}\,\,\|u(t)\|_{L^{q_{\mathrm{sc}}}}\leq \eta+1\right\},
\end{equation}
where $\delta>0$ and $\eta>0$ are constants to be chosen sufficiently small. We equip $\Xi$ with the metric
\begin{equation}\label{estiG}
    d_{\Xi}(u,v):=\sup_{t>0}t^\beta\|u(t)-v(t)\|_{L^q}+\sup_{t>0}\|u(t)-v(t)\|_{L^{q_{\mathrm{sc}}}},\qquad \hbox{for all}\,\, u,v\in\Xi.
\end{equation}
Then $(\Xi,d)$ is a nonempty complete metric space. For any $u\in\Xi$, define the mapping $\Phi$ by
\begin{equation}\label{estiH}
    \Phi(u)(t):=S_{\mathbb{H}}(t)u_0 + \int_{0}^{t}S_{\mathbb{H}}(t-s) \int_0^s(s-\tau)^{-\gamma}\abs{u(\tau)}^{p_1-1}u(\tau) \, \mathrm{d}\tau \,\mathrm{d}s + \int_{0}^{t}S_{\mathbb{H}}(t-s) |u(s)|^{p_2-1}u(s) \mathrm{d}s,
\end{equation}
for all $t\geq0$. Let us prove that  $\Phi: \Xi \rightarrow \Xi$. Using the assumptions $q\geq p_i$ for $i=1,2$, the estimate \eqref{estiE}, the definition of $\Xi$ in \eqref{estiF}, and Lemma \ref{Lp-Lqestimate}, we obtain for any $u \in \Xi$, 
\begin{eqnarray}\label{estiI}
    t^\beta\|\Phi(u)(t)\|_{L^q}&\leq& \eta+C\delta^{p_1} t^\beta\int_0^t\int_0^s(t-s)^{-\frac{Q(p_1-1)}{2
    q}}(s-\sigma)^{-\gamma}\sigma^{-\beta p_1}\,d\sigma\,ds\nonumber\\
    &{}&+\,C\delta^{p_2} t^\beta\int_0^t(t-s)^{-\frac{Q(p_2-1)}{2
    q}}s^{-\beta p_2}\,ds.
\end{eqnarray}
Next, by applying estimates \eqref{estiA} and \eqref{estiD}, and using the assumptions $p_2<p_1$ and $p_i\beta<1$ for $i=1,2$, we evaluate the integral terms appearing in \eqref{estiI}. Specifically, we obtain
\begin{equation}\label{estiJ}
 \int_0^t\int_0^s(t-s)^{-\frac{Q(p_1-1)}{2
    q}}(s-\sigma)^{-\gamma}\sigma^{-\beta p_1}\,d\sigma\,ds= C\int_0^t(t-s)^{-\frac{Q(p_1-1)}{\beta
    q}}s^{1-\gamma-\beta p_1}\,ds=C t^{-\beta},
\end{equation}
and
\begin{equation}\label{estiJ+}
 \int_0^t(t-s)^{-\frac{Q(p_2-1)}{2
    q}}s^{-\beta p_2}\,ds=Ct^{1-\frac{Q(p_2-1)}{2
    q}-\beta p_2}=C t^{-\beta},
\end{equation}
for all $t\geq0.$ Substituting \eqref{estiJ}-\eqref{estiJ+} into inequality \eqref{estiI}, we conclude
\begin{equation}\label{estiK}
    t^\beta\|\Phi(u)(t)\|_{L^q}\leq \eta+C\delta^{p_1}+C\delta^{p_2}.
\end{equation}
For the $L^{q_{\mathrm{sc}}}$ norm, again using Lemma \ref{Lp-Lqestimate}, estimate \eqref{estiEE}, and the assumptions $q\geq p_i$ for $i=1,2$, we have
$$
  \|\Phi(u)(t)\|_{L^{q_{\mathrm{sc}}}} \leq\eta+C\delta^{p_1} \int_0^t\int_0^s(t-s)^{-a_1}(s-\sigma)^{-\gamma}\sigma^{-\beta p_1}\,d\sigma\,ds+\,C\delta^{p_2} \int_0^t(t-s)^{-a_2}s^{-\beta p_2}\,ds,
$$
where
$$a_i:=\frac{Q}{2}\left(\frac{p_i}{q}-\frac{1}{q_{\mathrm{sc}}}\right)\qquad i=1,2.$$
Since $a_i<1$, $p_i\beta<1$ for $i=1,2$, and using the balance condition
$$2-a_1-\gamma-\beta p_1=1-a_2-\beta p_2=0,$$
we conclude,
\begin{equation}\label{estiAAA}
  \|\Phi(u)(t)\|_{L^{q_{\mathrm{sc}}}} \leq\eta+C\delta^{p_1} +\,C\delta^{p_2}.
\end{equation}
Combining estimates \eqref{estiK} and \eqref{estiAAA}, we see that if $\eta$ and $\delta$ are chosen sufficiently small so that 
$$\eta+C\delta^{p_1}+C\delta^{p_2}\leq\delta\quad\hbox{and}\quad C\delta^{p_1}+C\delta^{p_2}\leq 1,$$
then the mapping $\Phi$  is closed under $\Xi$, that is, $\Phi:\Xi\rightarrow\Xi$. Moreover, similar estimates show that $\Phi$ is a strict contraction on $\Xi$ (under the metric $d_\Xi$) when $\eta$ and $\delta$ are sufficiently small. Hence, by the Banach fixed-point theorem, there it exists a unique fixed point $u\in\Xi$, which is a mild solution of \eqref{1}.

We now establish that $u\in C([0,\infty),C_0(\mathbb{H}^{n}))$.  We begin by proving that $u\in C([0,T],C_0(\mathbb{H}^{n}))$ for some sufficiently small $T>0$. Indeed, the previous argument ensures uniqueness in $\Xi_T,$ where for any $T>0,$
$$
\Xi_T:=\left\{u\in
    L^\infty((0,T),L^q(\mathbb{H}^{n}));\;\sup_{0<t<T}t^\beta\|u(t)\|_{L^q}\leq\delta\right\}.
$$
Let $\tilde{u}$ denote the local mild solution of \eqref{1}
obtained in Theorem \ref{Local}. From inequality \eqref{estiB}, we have $q_{sc}<q<\infty$, which implies 
$$u_0\in C_0(\mathbb{H}^{n})\cap L^{q_{sc}}(\mathbb{H}^{n})\subset C_0(\mathbb{H}^{n})\cap L^{q}(\mathbb{H}^{n}).$$ 
Therefore,  Theorem \ref{Local} guarantees that $\tilde{u}\in C([0,T_{\max}),C_0(\mathbb{H}^{n})\cap L^q(\mathbb{H}^{n}))$. Consequently, we have
$$\sup\limits_{t\in(0,T)}t^\beta\|\tilde{u}(t)\|_{L^q}\leq\delta$$
for sufficiently small  $T>0$. By uniqueness, $u=\tilde{u}$ on
$[0,T]$, and hence $u\in C([0,T],C_0(\mathbb{H}^{n}))$.\\
Next, we extend the regularity to $[T,\infty)$ via a bootstrap argument. For $t>T,$ we express $u(t)$ as
$$
  u(t)-S_{\mathbb{H}}(t)u_0= I_1(t)+I_2(t)+I_3(t)+I_4(t),
$$
where
$$I_1(t):= \int_0^tS_{\mathbb{H}}(t-s)\int_0^T
    (s-\sigma)^{-\gamma}\abs{u}^{p_1-1}u(\sigma)\;d\sigma\,ds,\quad I_2(t):= \int_0^tS_{\mathbb{H}}(t-s)\int_T^s(s-\sigma)^{-\gamma}\abs{u}^{p_1-1}u(\sigma)\;d\sigma\,ds,$$
   $$I_3(t):= \int_0^TS_{\mathbb{H}}(t-s)\abs{u}^{p_2-1}u(s)\,ds\quad\hbox{and}\quad I_4(t):=\int_T^tS_{\mathbb{H}}(t-s)\abs{u}^{p_2-1}u(s)\,ds .$$
Since $u\in C([0,T],C_0(\mathbb{H}^{n})),$ it follows that $I_1,I_3\in
C([T,\infty),C_0(\mathbb{H}^{n})).$ Moreover, sing estimates similar to those in the fixed-point construction and noting $t^{-\beta}\leq
T^{-\beta}<\infty$, we also obtain $I_1,I_3\in C([T,\infty),L^q(\mathbb{H}^{n}))$. Now, from inequality \eqref{estiB}, we have $q>Q(p_1-1)/2$, which guarantees the existence of $r\in(q,\infty]$ satisfying
\begin{equation}\label{estiL}
\frac{Q}{2}\left(\frac{p_1}{q}-\frac{1}{r}\right)<1.
\end{equation}
For $T<s<t$, since $u\in
L^\infty((0,\infty),L^q(\mathbb{H}^{n})),$ we deduce
$$|u|^{p_1-1}u\in L^\infty((T,s),L^{q/p_1}(\mathbb{H}^{n})),\quad\quad |u|^{p_2-1}u\in L^\infty((T,t),L^{q/p_2}(\mathbb{H}^{n})).$$
 Applying (\ref{semigroup}), (\ref{estiL}), and the fact that $p_2<p_1$, it follows that $I_2, I_4\in
C([T,\infty),L^r(\mathbb{H}^{n}))$.
Since the terms $S_{\mathbb{H}}(t)u_0$, $I_1$, and $I_3$ belong to
$C([T,\infty),C_0(\mathbb{H}^{n}))\cap
C([T,\infty),L^q(\mathbb{H}^{n})),$ we conclude that $u\in
C([T,\infty),L^r(\mathbb{H}^{n}))$. We ultimately deduce, using a bootstrap argument, that $u\in
C([T,\infty),C_0(\mathbb{H}^{n})).$ \\
On the other hand, as $u\in C([0,T],C_0(\mathbb{H}^{n}))$ for all $T>0$, one can verify that the function
$$f(t):= \int_0^t(t-s)^{-\gamma}\abs{u(s)}^{p_1-1}u(s) \, \mathrm{d}s  +  |u(t)|^{p_2-1}u(t),\qquad \hbox{for all}\,\,t\in (0,T),$$
 belongs to $L^1((0,T),L^{q_{\mathrm{sc}}}(\mathbb{H}^{n}))$. By applying \cite[Lemma 4.1.5]{CH} and using the continuity of the semigroup $S_{\mathbb{H}}(t)$,  we deduce that $u\in C([0,T],L^{q_{\mathrm{sc}}}(\mathbb{H}^{n}))$ for any $T>0$, that is, $u\in C([0,\infty),L^{q_{\mathrm{sc}}}(\mathbb{H}^{n}))$.\\
\noindent{\it Step 2. Uniqueness.} Let $u,v\in C([0,\infty),C_0(\mathbb{H}^{n})\cap L^{q_{\mathrm{sc}}}(\mathbb{H}^{n}))$ be two global mild solutions. Then for any finite $T>0$, we have $u,v\in C([0,T],C_0(\mathbb{H}^{n})\cap L^{q_{\mathrm{sc}}}(\mathbb{H}^{n}))$. Fix $T>0$, and let $r\in\{q_{\mathrm{sc}},\infty\}$. We consider the difference $u(t)-v(t)$, which satisfies the integral inequality
\begin{eqnarray*}
  \|u(t)-v(t)\|_{L^r}&\leq& C(p,\gamma)
  \int_0^t\int_0^s(s-\sigma)^{-\gamma}
  \|u(\sigma)-v(\sigma)\|_{L^r}\,d\sigma\,ds+C(p)
  \int_0^t
  \|u(s-v(s)\|_{L^{r}}\,ds\\
   &=&C(p_1,\gamma)
   \int_0^t\int_\sigma^t(s-\sigma)^{-\gamma}\|u(\sigma)-
   v(\sigma)\|_{L^{r}}\,ds\,d\sigma +C(p_2)
  \int_0^t
  \|u(s-v(s)\|_{L^{r}}\,ds\\
   &=&C(p_1,\gamma)\int_0^t(t-\sigma)^{1-\gamma}\|u(\sigma)-
   v(\sigma)\|_{L^{r}}\,d\sigma+C(p_2)
  \int_0^t
  \|u(s-v(s)\|_{L^{r}}\,ds\\
  &\leq&\left(C(p_1,\gamma)T^{1-\gamma}+C(p_2)
\right)  \int_0^t
  \|u(s-v(s)\|_{L^{r}}\,ds.
\end{eqnarray*}
Applying Gronwall's inequality, we conclude that $u(t)=v(t)$ for every $t\in[0,T]$. Since $T>0$ was arbitrary, the uniqueness holds globally. This completes the proof in the case $p_2= \tilde{p}_2$.\\							
\noindent {\it The Case $p_2> \tilde{p}_2$.}  We treat this case using the comparison principle established in Theorem \ref{comparison}. We aim to show that the mild solution $u$ of \eqref{1} satisfies $0\leq u(t)\leq v(t)$, for all $t\geq0$, where $v$ denotes the mild solution of the same equation with the exponent $p_2$ replaced by $\tilde{p}_2$. From the case $p_2=\tilde{p}_2$, the solution $v$ exists globally and is locally bounded. Thus, by Theorem \ref{Local} and the comparison principle $u$ also exists globally under suitable smallness assumptions on the initial data. Indeed, let the initial data $u_0\geq0$ satisfy $0\leq u_0<\epsilon_0$ for some sufficiently small constant $\epsilon_0\ll1$. Define the associated mild solution operators $\Phi$ and $\tilde{\Phi}$ corresponding to the exponents $p_2$ and $\tilde{p}_2$, respectively, as
 		$$\Phi(u)(t):=S_{\mathbb{H}}(t)u_0 + \int_{0}^{t}S_{\mathbb{H}}(t-s) \int_0^s(s-\tau)^{-\gamma}(u(\tau))^{p_1}\, \mathrm{d}\tau \,\mathrm{d}s + \int_{0}^{t}S_{\mathbb{H}}(t-s)(u(s))^{p_2} \mathrm{d}s,$$
	$$\tilde{\Phi}(v)(t):=S_{\mathbb{H}}(t)u_0 + \int_{0}^{t}S_{\mathbb{H}}(t-s) \int_0^s(s-\tau)^{-\gamma}(v(\tau))^{p_1}\, \mathrm{d}\tau \,\mathrm{d}s + \int_{0}^{t}S_{\mathbb{H}}(t-s)(v(s))^{\tilde{p}_2} \mathrm{d}s.$$
	In fact, the solutions $u$ and $v$ are obtained as the limits of the following Picard sequences
	\begin{equation} \label{pic1A}
		u_{1}(t) =S_\mathbb{H}(t)u_0 ,\qquad u_{k}(t) = \Phi (u_{k-1})(t),\,\,\,k\geq2 , 
	\end{equation}
	\begin{equation} \label{pic12A}
		v_{1}(t) =S_\mathbb{H}(t)u_0 ,\qquad	 v_{k}(t) = \tilde{\Phi} (v_{k-1})(t),\,\,\, k\geq2 .
	\end{equation}
 We observe that both solutions begin with the same initial data, that is, $0\leq u_1= v_1$, and $v_1(t) \leq \|u_0\|_{L^\infty}<\epsilon_0$. Because $\tilde{p}_2<p_2$, we deduce
		$$u_1^{p_2}=v_1^{p_2}\leq  v_1^{\tilde{p}_2}.$$ 
		Hence, by Theorem \ref{comparison} (comparison principle), it follows that $0\leq u_2\leq v_2$.\\ 
		On the other hand, as established in the previous case, the solution $v$ is global, and its iterates $(v_n)_n\subset\Xi$, where $\Xi$ is the function space defined in \eqref{estiF}. Recalling that $\tilde{p}_2=(p_1+1-\gamma)/(2-\gamma)$,  we estimate $ \|v_2(t)\|_{L^{\infty}}$ as follows
		\begin{eqnarray*}
		 \|v_2(t)\|_{L^{\infty}}&\leq & C t^{-\frac{Q}{2q_{\mathrm{sc}}}} \|u_0\|_{L^{q_{\mathrm{sc}}}} +
		C\int_{0}^{t} (t-s)^{-\frac{Q}{2q}p_1} \int_0^s (s-\tau)^{-\gamma} \|v_1(\tau)\|^{p_1}_{L^{q}} \, \mathrm{d}\tau \,\mathrm{d}s  \\
		&{}& +\, C \int_{0}^{t}(t-s)^{-\frac{Q}{2q}\tilde{p}_2} \|v_1(s)\|^{\tilde{p}_2}_{L^{q}} \mathrm{d}s \\
		&\leq& C t^{-\frac{Q}{2q_{\mathrm{sc}}}} \|u_0\|_{L^{q_{\mathrm{sc}}}}  + C \delta^{p_1}  t^{2-\frac{Q}{2q}p_1 - \gamma - \beta p_1} + C \delta^{\tilde{p}_2}  t^{1-\frac{Q}{2q} \tilde{p}_2-\beta \tilde{p}_2} .
	\end{eqnarray*}
From the choice of the parameters $q_{\mathrm{sc}}$, $\beta$, and the smallness of $\delta$ and $\|u_0\|_{L^{q_{\mathrm{sc}}}}$, one can deduce the following bound for large times $t\geq1$,
		\begin{equation}\label{vn1}
			\sup_{t\geq1}\|v_2(t)\|_{L^{\infty}}  \leq \epsilon_0+\delta.
		\end{equation}
		For small times $t\in(0,1]$, since $\|v_1(t)\|_{L^{\infty}}<\epsilon_0$, we obtain
		\begin{equation*}
			\|v_2(t)\|_{L^{\infty}}  \leq  \|u_0\|_{L^{\infty}}  +  \frac{\epsilon_0^{p_1}}{(1-\gamma)(2-\gamma)}  + \epsilon_0^{\tilde{p}_2}  ,
		\end{equation*}
		so that, for $\epsilon_0$ sufficiently small,
		\begin{equation}\label{vn2}
		\sup_{0<t\leq1}\|v_2(t)\|_{L^{\infty}}  \leq 2\epsilon_0.
			\end{equation}
		Combining estimates \eqref{vn1} and \eqref{vn2}, and choosing $\delta,\epsilon_0$ small enough, we infer that
	$$
		\|v_{2}(t)\|_{L^{\infty}}  \leq \max\{\epsilon_0+\delta,2\epsilon_0\}<1.
	$$
	Hence, 
		$$u_2^{p_2}(t)\leq v_2^{p_2}(t)\leq  v_2^{\tilde{p}_2}(t),\quad\hbox{for all}\,\,t>0.$$ 
	Applying Theorem \ref{comparison} again, we conclude $0\leq u_3(t)\leq v_3(t)$. By induction and appropriate choice of small constants $\delta,\epsilon_0$, we obtain 
	$$0\leq u_n(t)\leq v_n(t)<1,\quad\hbox{for all}\,\,n\in\N,\,t>0,$$
	which implies $\|u(t)\|_{L^{\infty}}\leq  \|v(t)\|_{L^{\infty}}\leq 1$, for all $t>0$. Thus, $\|u\|_{L^\infty((0,t)\times\mathbb{H}^{n})}<\infty$ for all $t>0$ and the global existence of $u$ follows from Theorem \ref{Local}.\\	
	
\noindent {\it The Case $p_2^{**}<p_2< \tilde{p}_2$.} Since the inequality $p_2<\tilde{p}_2$ is equivalent to $p_1 > (p_2-1)(2-\gamma)+1 =:\tilde{p}_1$, it follows that $w^{p_1}\leq w^{\tilde{p}_1}$ whenever $w<\varepsilon_0$. Moreover, we observe that $(\tilde{p}_1+1-\gamma)/(2-\gamma) = p_2$, which allows us to  apply the same reasoning as in the previous case, by interchanging the roles of $p_1$ and $p_2$, using $\tilde{p}_1$ in place of $\tilde{p}_2$. That is, we obtain the estimate $\|u(t)\|_{L^{\infty}}\leq  \|v(t)\|_{L^{\infty}}$, for all $t>0$, where $v$ is the solution of the integral equation
$$v(t)=S_{\mathbb{H}}(t)u_0 + \int_{0}^{t}S_{\mathbb{H}}(t-s) \int_0^s(s-\tau)^{-\gamma}(v(\tau))^{\tilde{p}_1}\, \mathrm{d}\tau \,\mathrm{d}s + \int_{0}^{t}S_{\mathbb{H}}(t-s)(v(s))^{p_2} \mathrm{d}s.$$
Noting that the assumption $p_2^{**}<p_2$ is equivalent to $\tilde{p}_1>p_1^*$, we conclude that the global existence of $u$ holds in this case as well.

\section{Finite-time blow up}\label{sec4}
This section is dedicated to the proof of Theorem \ref{T4}. Let $v$ be the mild solution of the following nonlocal problem:
\begin{equation}\label{11}
	\left\{\begin{array}{ll}
		\displaystyle v_{t}-\Delta_{\mathbb{H}} v =  \int_0^t(t-s)^{-\gamma}v^{p_1}(s)\,\mathrm{d}s, &\displaystyle {\eta\in {\mathbb{H}}^{n},\,t>0,}\\
		{}\\
		\displaystyle{v(\eta,0)=  u_0(\eta),\qquad\qquad}&\displaystyle{\eta\in {\mathbb{H}}^{n},}
	\end{array}
	\right.
\end{equation}
and let $w$ be the mild solution of the local problem:
\begin{equation}\label{111}
	\left\{\begin{array}{ll}
		w_{t}-\Delta_{\mathbb{H}} w =  w^{p_2}, &\displaystyle {\eta\in {\mathbb{H}}^{n},\,t>0,}\\
		{}\\
		\displaystyle{w(\eta,0)=  u_0(\eta),\qquad\qquad}&\displaystyle{\eta\in {\mathbb{H}}^{n}.}
	\end{array}
	\right.
\end{equation}
Using the nonnegativity of the solution $u\geq0$, it follows that $u$ is a mild supersolution to both problems \eqref{11} and \eqref{111}. Applying the comparison principle (Theorem \ref{comparison}), we conclude that  $u\geq v$ and $u\geq w$.\\
If $p_1\leq p_1^*$, then by \cite[Theorem~2]{FinoKiraneTorebek} and the fact that every mild solution is also a weak solution, the function $v$ blows up in finite time, and thus $u$ must also blow up. Similarly, if $p_2\leq p_2^*$, the result of \cite[Theorem~3.1]{Pohozaev} implies that $w$ blows up in finite time, leading to blow-up of $u$ as well. This concludes the proof of Theorem \ref{T4}.


\section*{Conclusion}

This article deals with a semilinear heat equation on the Heisenberg group, incorporating a mixed local and nonlocal nonlinearity. The main results concern the local and global solvability of mild solutions for regular, nonnegative initial data. Using the Banach fixed point theorem combined with appropriate $L^p-L^q$ estimates, we proved the existence and uniqueness of local-in-time solutions. Furthermore, under suitable growth conditions on the nonlinear terms, we extended the existence to global-in-time solutions.

To complement these results, we employed the capacity method to demonstrate the nonexistence of global solutions in the supercritical regime. In addition, sharp lifespan estimates were obtained in the supercritical regime, providing insight into how the initial data influences the time of blow-up.

The results obtained in this work extend existing literature by incorporating mixed local-nonlocal nonlinearities within the geometric framework of the Heisenberg group. In contrast to the classical Euclidean setting, the use of sub-Riemannian structure and stratified group properties introduces additional challenges and rich dynamics. The use of fixed point arguments in this context highlights the effectiveness of abstract functional methods in non-Euclidean analysis.

The novelty of this work lies in formulating and proving Fujita-type results for equations involving both memory-type and instantaneous nonlinearities on the Heisenberg group. These findings contribute to the ongoing effort to understand how geometry and nonlocal interactions influence solution behavior in nonlinear parabolic problems.

Practically, the results provide a theoretical foundation for analyzing more complex models that appear in physics and geometry, where nonlocal effects and memory terms are relevant. Future research may explore extensions involving time-space fractional operators, more singular initial data, or the influence of different types of nonlocal terms on other classes of Lie groups or manifolds.

\section*{Acknowledgments}

Ahmad Z. Fino is supported by the Research Group Unit, College of Engineering and Technology, American University of the Middle East. 

\section*{Author Contributions}

A.Z.~Fino carried out the main analysis and computations of the manuscript. M.~Kirane contributed through discussions on the formulation of the problem, provided key ideas, and assisted with references and manuscript revision. All authors contributed equally to this work. All authors have read and approved the final version of the manuscript.
\section*{Conflict of Interest}
The authors declare no conflict of interest.

\par\bigskip


\bibliographystyle{elsarticle-num}


\bigskip
\begin{center}
\emph{Author Information}
\end{center}

\textbf{Ahmad Z. Fino} --- Doctor of applied mathematics, Associate Professor, College of Engineering and Technology, American University of the Middle East, Egaila 54200, Kuwait; e-mail: \emph{ahmad.fino@aum.edu.kw}; https://orcid.org/0000-0001-7067-1553. \vspace{6pt}

\textbf{Mokhtar Kirane}  --- Doctor of applied mathematics, Professor, Department of Mathematics, College of Computing and Mathematical Sciences, Khalifa University, P.O. Box: 127788,  Abu Dhabi, UAE; e-mail: \emph{mokhtar.kirane@ku.ac.ae}; https://orcid.org/0000-0002-4867-7542.\vspace{6pt}

\textbf{Zineb Sabbagh} --- Doctor of applied mathematics, Department of Mathematics, University M'hamed Bougara Boumerdes, Algeria; e-mail: \emph{szinebedp@gmail.com}

\end{document}